\documentclass[12pt, reqno]{amsart}
\usepackage{amsmath, amsthm, amscd, amsfonts, amssymb}
\usepackage[bookmarksnumbered]{hyperref}

\newtheorem{theorem}{Theorem}[section]

\theoremstyle{definition}
\newtheorem{definition}[theorem]{Definition}
\newtheorem{example}[theorem]{Example}

\theoremstyle{remark}

\numberwithin{equation}{section}

\begin{document}

\title[]{Fuzzy stability of additive--quadratic functional equations }
\author[M. Eshaghi Gordji,  N.Ghobadipour and  J. M. Rassias]{M. Eshaghi Gordji,  N.Ghobadipour and J. M. Rassias}
\address{Department of Mathematics, Semnan University, P. O. Box 35195-363, Semnan, Iran and Section of Mathematics and Informatics,
Pedagogical Department, National and Capodistrian University of
Athens, 4, Agamemnonos St., Aghia Paraskevi, Athens 15342, Greece
} \email{madjid.eshaghi@gmail.com,    ghobadipour.n@gmail.com,
jrassias@primedu.uoa.gr}

 \subjclass[2000]{Primary 46S40; Secondary 39B52, 39B82, 26E50,
 46S50.}

\keywords{Fuzzy normed space; additive functional equation;
quadratic functional equation; fuzzy generalized Hyers-Ulam
stability.}

\begin{abstract} In this paper we investigate the generalized
Hyers- Ulam stability of the functional equation $$f (2x +y)+f (2x
-y)= f (x + y)+ f (x -y)+2f (2x)-2f (x)$$ in fuzzy Banach spaces.
\end{abstract}
\maketitle

%--------------------------------------------------------------------------------------%

\section{Introduction and preliminaries }
In 1984, Katsaras \cite{Kat} defined a fuzzy norm on a linear
space and at the same year Wu and Fang \cite{Co} also introduced a
notion of fuzzy normed space and gave the generalization of the
Kolmogoroff normalized theorem for  fuzzy topological linear
space. In \cite{Bi}, Biswas defined and studied fuzzy inner
product spaces in  linear space. Since then some mathematicians
have defined fuzzy metrics and norms on a linear space from
various points of view \cite{Ba,Fe,Kri,Shi,Xi}. In 1994, Cheng and
Mordeson introduced a definition of fuzzy norm on a linear space
in such a manner that the corresponding induced fuzzy metric is of
Kramosil and Michalek type \cite{Kra}. In 2003, Bag and Samanta
\cite{Ba} modified the definition of Cheng and Mordeson \cite{Che}
by removing a regular condition. They also established a
decomposition theorem of a fuzzy norm into a family of crisp norms
and investigated some properties of fuzzy norms (see \cite{Bag2}).
Following \cite{Bag1}, we give the employing notion of a fuzzy
norm.\\
Let X be a real linear space. A function $N : X \times \Bbb R
\longrightarrow [0,1]$(the so-called fuzzy subset) is said to be a
fuzzy norm on X if for all $x, y \in X$ and all $a,b \in \Bbb
R$:\\
$(N_1)~~N(x,a)=0$ for $a\leq 0;$\\
$(N_2)~~x=0$ if and only if $N(x,a)=1$ for all $a>0;$\\
$(N_3)~~N(ax,b)=N(x,\frac{b}{|a|})$ if $a\neq0$;\\
$(N_4)~~N(x+y,a+b)\geq min \{N(x,a), N(y,b)\};$\\
$(N_5)~~N(x,.)$ is non-decreasing function on $\Bbb R$ and $\lim
_{a \to \infty} N(x,a)=1;$\\
$(N_6)~~$ For $x\neq 0,$ $N(x,.)$ is (upper semi) continuous on
$\Bbb R.$\\
The pair $(X,N)$ is called a fuzzy normed linear space. One may
regard $N(x,a)$ as the truth value of the statement 'the norm of
$x$
is less than or equal to the real number $a^{'}$.\\
\begin{example}
Let $(X,\|.\|)$ be a normed linear space. Then $$N(x,a)=\left\{%
\begin{array}{ll}
    \frac{a}{a+\|x\|}, & a>0~~~, x \in X, \\
    0,~~~ & a\leq0, x \in X \\
\end{array}%
\right.    $$

is a fuzzy norm on X.
\end{example}

\begin{definition}
Let $(X,N)$ be a fuzzy normed linear space. Let ${x_n}$ be a
sequence in X. Then ${x_n}$ is said to be convergent if there
exists $x \in X$ such that $\lim_{n \to \infty}N(x_n-x,a)=1$ for
all $a>0.$ In that case, $x$ is called the limit of the sequence
${x_n}$ and we denote it by $N - lim_{n \to \infty} ~ x_n = x.$
\end{definition}
\begin{definition}
A sequence ${x_n}$ in X is called Cauchy if for each $\epsilon >
0$ and each $a
> 0$ there exists $n_0$ such that for all $n\geq n_0$ and all $p > 0,$ we have
$N(x_{n+p} - x_n, a) > 1$ - $\epsilon$.
\end{definition}
It is known that every convergent sequence in  fuzzy normed space
is Cauchy. If each Cauchy sequence is convergent, then the fuzzy
norm is said to be complete and the fuzzy normed space is called a
fuzzy Banach space.\\
The study of the stability problem of functional equations was
introduced by Ulam  \cite {Ul}  Let $(G_1,.)$ be a group and let
$(G_2,*)$ be a metric group with the metric $d(.,.).$ Given
$\epsilon >0$, does there exist a $\delta
>0$, such that if a mapping $h:G_1\longrightarrow G_2$ satisfies the
inequality $d(h(x.y),h(x)*h(y)) <\delta$ for all $x,y\in G_1$,
then there exists a homomorphism $H:G_1\longrightarrow G_2$ with
$d(h(x),H(x))<\epsilon$ for all $x\in G_1?$ In other words, under
what condition does there exist a homomorphism near an approximate
homomorphism? The concept of stability for functional equations
arises when we replace the functional equation by an inequality
which acts as a perturbation of the equation. In 1941, D. H. Hyers
\cite{Hy} gave the first affirmative  answer to the question of
Ulam for Banach spaces $E$ and $E'$. Let
$f:{E}\longrightarrow{E'}$ be a mapping between Banach spaces such
that
$$\|f(x+y)-f(x)-f(y)\|\leq \delta $$
for all $x,y\in E,$ and for some $\delta>0.$ Then there exists a
unique additive mapping $T:{E}\longrightarrow{E'}$ such that
$$\|f(x)-T(x)\|\leq \delta$$
for all $x\in E.$ Now assume that $E$ and $E^{'}$ are real normed
spaces with $E^{'}$ complete, $f:E \to E^{'}$ is a mapping such
that  the mapping $t\rightarrowtail f(tx)$ is continuous in  $t
\in \Bbb R$ for each fixed $x \in E,$ and that there exist
$\delta\geq 0$ and $p\neq 1$ such that
$$\|f(x+y)-f(x)-f(y)\|\leq \delta (\|x\|^p+\|y\|^p)$$ for all $x,y
\in E.$ Then there exists a unique linear map $T:E \to E^{'}$ such
that
$$\|f(x)-T(x)\| \leq \frac{2 \delta \|x\|^p}{|2^p-2|}$$ for all $x
\in E.$  (see \cite{TRa}).\\
On the other hand J. M. Rassias \cite {JRa1,JRa2,JRa3,JRa4}
generalized the Hyers stability result by presenting a weaker
condition controlled by a product of different powers of norms.
According to J. M. Rassias Theorem:
\begin{theorem}\label{t2}
If it is assumed that there exist constants $\Theta\geq0$ and
$p_1, p_2 \in \Bbb R$ such that $p=p_1+p_2\neq 1,$ and $f:E\to
E^{'}$ is a map from a norm space $E$ into a Banach space $E^{'}$
such that the inequality $$\|f(x+y)-f(x)-f(y)\|\leq \epsilon
\|x\|^{p_1} \|y\|^{p_2}$$ for all $x,y\in E,$ then there exists a
unique additive mapping $T:E\to E^{'}$ such that $$\|f(x)-T(x)\|
\leq \frac{\Theta}{2-2^p}\|x\|^p ,$$ for all $x\in E.$ If in
addition for every $x\in E,$ $f(tx)$ is continuous in  $t \in \Bbb
R$ for each fixed $x,$ then $T$ is linear.
\end{theorem}
Quadratic functional equation was used to characterize inner
product spaces \cite{Ac,Am,Jo}. Several other functional equations
were also used to characterize inner product spaces. A square norm
on an inner product space satisfies the important parallelogram
equality\\
$$\|x+y\|^2+\|x-y\|^2=2(\|x\|^2+\|y\|^2).$$
The functional equation $$f(x+y)+f(x-y)=2f(x)+2f(y) \eqno(1.1)$$
is related to a symmetric bi-additive function \cite{Ac,Kan}. It
is natural that each equation is called a quadratic functional
equation. In particular, every solution of the quadratic equation
(1.1) is said to be a quadratic function. It is well known that a
function f between real vector spaces is quadratic if and only if
there exists a unique symmetric bi-additive function B such that
$f(x)=B(x,x)$ for all $x$ (see \cite{Ac,Kan}). The bi-additive
function B is given by
$$B(x,y)=\frac{1}{4}(f(x+y)-f(x-y)).\eqno(1.2)$$ A Hyers-Ulam
stability problem for the quadratic functional equation (1.1) was
proved by Skof for functions $f:E_1 \longrightarrow¨E_2$ where
$E_1$ is a normed space and $E_2$ is a Banach space (see
\cite{Sk}). Cholewa \cite{cho} noticed that the theorem of Skof is
still true if the relevant domain $E_1$ is replaced by an Abelian
group. In the paper \cite{Cz}, Czerwik proved the generalized
Hyers-Ulam stability of the quadratic functional equation (1.1).
Grabiec \cite{Gr} has generalized these results mentioned above.
Jun and Lee \cite{Ju} proved the generalized Hyers-Ulam stability
of the
pexiderized quadratic equation (1.1).\\
A. Najati and M.B. Moghimi \cite{Na}, have obtained the
generalized Hyers- Ulam stability for a functional equation
deriving from
quadratic and additive functions in quasi-Banach spaces.\\
In this paper, we deal with the the following functional equation
deriving from quadratic and additive functions:
$$f (2x +y)+f (2x -y)= f (x + y)+ f (x -y)+2f (2x)-2f (x).
\eqno(1.3)$$  It is easy to see that the function $f (x)= ax^2 +
bx+ c$ is a solution of the functional equation $(1.3).$ The main
purpose of this paper is to establish some versions of the
generalized Hyers- Ulam stability for the function  equation
$(1.3)$ in  fuzzy normed linear spaces.

\section{Main result}
Throughout this section, assume that $X,$ $(Y,N)$ and $(Z,N^{'})$
are linear space, fuzzy normed space and fuzzy Banach space,
respectively. We start our works with a fuzzy generalized Hyers-
Ulam type theorem for the functional equation $(1.3).$
\begin{theorem}\label{t2}
Let $\varphi_1: X\times X \to Z$ be a function such that for some
$0 < \alpha < 4$
$$N^{'}(\varphi_1(\frac{2x}{3},2y),a)\geq N^{'}(\alpha
\varphi_1(\frac{x}{3},y),a) \eqno(2.1)$$ for all $x \in X,$ $y \in
\{0,\frac{x}{3},\frac{4x}{3},\frac{-2x}{3},x\}$ and $a>0,$ and
$\lim_{n \to \infty}N^{'}(\varphi_1(2^nx,2^ny),4^na)=1$ for all
$x, y \in X$ and $a
> 0.$ Let $f: X \to Y$ be
an even function with $f(0)=0$ satisfying
$$N(f(2x+y)+f(2x-y)-f(x+y)-f(x-y)-2f(2x)+2f(x),a)\geq N^{'}(\varphi_1(x,y),a) \eqno(2.2)$$ for all
$a
> 0$ and all $x, y \in X.$ Then there exists a unique quadratic mapping $Q:
X \to Y$ such that $$N(Q(x)-f(x),a)\geq
N^{''}_1(x,\frac{a(4-\alpha)}{6})\eqno(2.3)$$ for all $x \in X$
and all $a
> 0,$ where
\begin{align*}
N^{''}_1(x,a):=min \{&N^{'}(\varphi_1(\frac{x}{3},\frac{x}{3}),a)
,N^{'}(\varphi_1(\frac{x}{3},x),a),\\
&N^{'}(\varphi_1(\frac{x}{3},\frac{4x}{3}),a)
,N^{'}(\varphi_1(\frac{x}{3},\frac{-2x}{3}),a)
,N^{'}(\varphi_1(\frac{x}{3},0),a)\}.
\end{align*}
\end{theorem}
\begin{proof}
By replacing $y$ by $x+y$ in $(2.2),$ we get
$$N(f(3x+y)+f(x-y)-f(2x+y)-f(y)-2f(2x)+2f(x),a)\geq
N^{'}(\varphi_1(x,x+y),a)\eqno(2.4)$$ for all $x,y \in X$ and
$a>0.$ Replacing $y$ by $-y$ in $(2.4),$ we get
\begin{align*}
N(f(3x-y)&+f(x+y)-f(2x-y)-f(y)-2f(2x)+2f(x),a)\\
&\geq N^{'}(\varphi_1(x,x-y),a) \hspace{9.cm}(2.5)
\end{align*}
for all $x,y \in X$ and $a>0.$ It follows from $(2.2),$ $(2.4)$
and $(2.5),$
\begin{align*}
N(f(3x+y)&+f(3x-y)-2f(y)-6f(2x)+6f(x),3a)\\
&\geq min \{N^{'}(\varphi_1(x,y),a)
,N^{'}(\varphi_1(x,x+y),a),N^{'}(\varphi_1(x,x-y),a)\}
\hspace{2.cm}(2.6)
\end{align*}
for all $x,y \in X$ and $a>0.$ Letting $y=0$ in $(2.6),$ we  get
inequality
$$N(2f(3x)-6f(2x)+6f(x),3a)\geq min \{
N^{'}(\varphi_1(x,x),a),N^{'}(\varphi_1(x,0),a)\}\eqno(2.7)$$ for
all $x,y \in X$ and $a>0.$  Putting $y=3x$ in $(2.6),$ we get
\begin{align*}
N(f(6x)&-2f(3x)-6f(2x)+6f(x),3a)\\
&\geq min \{ N^{'}(\varphi_1(x,3x),a)
,N^{'}(\varphi_1(x,4x),a),N^{'}(\varphi_1(x,-2x),a)\}
\hspace{2.95cm}(2.8)
\end{align*}
$x,y \in X$ and $a>0.$ It follows from $(2.7)$ and $(N_3),$
$$N(-2f(3x)+6f(2x)-6f(x),3a)\geq min \{
N^{'}(\varphi_1(x,x),a),N^{'}(\varphi_1(x,0),a)\}\eqno(2.9)$$ for
all $x,y \in X$ and $a>0.$ Therefore we to obtain from $(2.8)$ and
$(2.9)$ the inequality
\begin{align*}
N(f(6x)-4f(3x),6a)\geq min \{
&N^{'}(\varphi_1(x,x),a),N^{'}(\varphi_1(x,3x),a),\\
&N^{'}(\varphi_1(x,4x),a),N^{'}(\varphi_1(x,-2x),a),N^{'}(\varphi_1(x,0),a)\}
\hspace{.75cm}(2.10)
\end{align*}
for all $x,y \in X$ and $a>0.$ If we replace $x$ by $\frac{x}{3}$
in $(2.10)$ for all $x \in X$ and $a>0,$ then we get then $$
N(f(2x)-4f(x),6a)\geq N^{''}_1(x,a)\eqno(2.11)$$ for all $x \in X$
and $a>0.$ Thus $$N(\frac{f(2x)}{4}-f(x),\frac{3a}{2})\geq
N^{''}_1(x,a)\eqno(2.12)$$ for all $x \in X$ and $a>0.$ Replacing
$x$ by $2^nx$ in $(2.12),$ we get
$$N(\frac{f(2^{n+1}x)}{4}-f(2^nx),\frac{3a}{2})\geq N^{''}_1(2^nx,a)\eqno(2.13)$$ for all $x\in
X$ and $a>0.$ Using $(2.1)$ we get
$$N(\frac{f(2^{n+1}x)}{4}-f(2^nx),\frac{3a}{2})\geq N^{''}_1(x,\frac{a}{\alpha^n})\eqno(2.14)$$
for all $x \in X$ and $a>0.$  Replacing $a$ by $\alpha^na$ we see
that
$$N(\frac{f(2^{n+1}x)}{4^{n+1}}-\frac{f(2^nx)}{4^n},\frac{3a\alpha ^n}{2(4^n)})
\geq N^{''}_1(x,a)\eqno(2.15)$$ for all $x \in X$ and $a>0.$ It
follows from
$\frac{f(2^nx)}{4^n}-f(x)=\sum_{i=0}^{n-1}\frac{f(2^{i+1}x)}{4^{i+1}}-\frac{f(2^ix)}{4^i}$
and $(2.15)$ that
$$N(\frac{f(2^nx)}{4^n}-f(x),\sum_{i=0}^{n-1}\frac{3a\alpha^i}{2(4^i)})
\geq min
\bigcup_{i=0}^{n-1}\{N(\frac{f(2^{i+1}x)}{4^{i+1}}-\frac{f(2^ix)}{4^i},\frac{3a\alpha^i}{2(4^i)})\}\\
\geq N^{''}_1(x,a)\eqno(2.16)$$ for all $x \in X$  and $a>0.$
Replacing $x$ with $2^mx$ in $(2.16)$ we observe that
$$N(\frac{f(2^{n+m}x)}{4^{n+m}}-\frac{f(2^mx)}{4^m},\sum_{i=0}^{n-1}\frac{3a\alpha
^i}{2(4^{i+m})})\geq N^{''}_1(2^mx,a)\geq
N^{''}_1(x,\frac{a}{\alpha^m}),$$ whence
$$N(\frac{f(2^{n+m}x)}{4^{n+m}}-\frac{f(2^mx)}{4^m},\sum_{i=m}^{n+m-1}\frac{3a\alpha
^i}{2(4^{i})})\geq N^{''}_1(x,a)$$ for all $x
\in X,$ $a>0$ and $m,n\geq 0.$\\
Hence
$$N(\frac{f(2^{n+m}x)}{4^{n+m}}-\frac{f(2^mx)}{4^m},a))
\geq N^{''}_1(x,\frac{a}{\sum_{i=m}^{n+m-1}\frac{3\alpha
^i}{2(4^{i})}})\eqno(2.17)$$ for all $x \in X,$ $a>0$ and $m,n\geq
0.$ Since $0<\alpha <4$ and
$\sum_{i=0}^{\infty}(\frac{\alpha}{4})^i< \infty$ the Cauchy
criterion for convergence and $(N_5)$ imply that
$\{\frac{f(2^nx)}{4^n}\}$ is a Cauchy sequence in $(Y,N).$ Since
$(Y,N)$ is a fuzzy Banach space, this sequence converges to some
point $Q(x)\in Y.$ So one can define the mapping $Q:X \to Y$ by
$Q(x):=N-\lim_{n \to \infty}\frac{f(2^nx)}{4^n}$ for all $x \in
X.$ Letting $m=0$ in $(2.17),$ we get
$$N(\frac{f(2^{n}x)}{4^{n}}-f(x),a)\geq N^{''}_1(x,\frac{a}
{\sum_{i=0}^{n-1}\frac{3\alpha ^i}{2(4^{i})}})\eqno(2.18)$$ for
all $x \in X$  and $a>0.$ Taking the limit as $n \to \infty$ and
using $(N_6)$ we get
$$N(Q(x)-f(x),a)\geq N^{''}_1(x,\frac{a(4-\alpha)}{6})$$
 for all $x \in X$ and $a>0.$ Now, we claim  that $Q$ is quadratic. Replace $x,y$ by
$2^nx,2^ny,$ respectively in $(2.2)$ to get
\begin{align*}
N(\frac{f(2^n(2x+y))}{4^n}&+\frac{f(2^n(2x-y))}{4^n}-\frac{f(2^n(x+y))}{4^n}
\\
&-\frac{f(2^n(x-y))}{4^n}-\frac{2f(2^n(2x))}{4^n}
+\frac{2f(2^nx)}{4^n},a)\\
&\geq N^{'}(\varphi_1(2^nx,2^ny),4^na)
\end{align*}
for all $x,y \in X$ and $a>0.$ Since $\lim_{n \to
\infty}N^{'}(\varphi_1(2^nx,2^ny),4^na)=1$ and $Q(0)=0,$ then by
Lemma $2.1$ of \cite{Na} we get that the mapping $Q: X \to Y$ is
quadratic.\\
To prove the uniqueness of $Q,$ let $Q^{'}:X \to Y$ be another
quadratic mapping satisfying $(2.3).$ Fix $x \in X.$ Clearly
$Q(2^nx)=4^nQ(x)$ and $Q^{'}(2^nx)=4^nQ^{'}(x)$ for all $n \in
\Bbb N.$ It follows from $(2.3)$ that
\begin{align*}
N(Q(x)-Q^{'}(x),a)&=N(\frac{Q(2^nx)}{4^n}-\frac{Q^{'}(2^nx)}{4^n},a)\\
&\geq min \{N(\frac{Q(2^nx)}{4^n}-\frac
{f(2^nx)}{4^n},\frac{a}{2}),
N(\frac{f(2^nx)}{4^n}-\frac{Q^{'}(2^nx)}{4^n},\frac{a}{2})\}\\
&\geq N^{''}_1(2^nx,\frac{a(4-\alpha)(4^n)}{12}) \geq
N^{''}_1(x,\frac{a(4-\alpha)(4^n)}{12\alpha^n})
\end{align*}
for all $x \in X$  and $a>0.$\\
Since $\lim_{n \to
\infty}\frac{a(4-\alpha)(4^n)}{12\alpha^n}=\infty,$ we obtain
$\lim_{n \to
\infty}N^{''}_1(x,\frac{a(4-\alpha)(4^n)}{12\alpha^n})=1.$
Therefore, $N(Q(x)-Q^{'}(x),a)=1$ for all $x \in X$ and all $a>0,$
whence $Q(x)=Q^{'}(x).$
\end{proof}
\begin{theorem}\label{t2}
Let $\varphi_2: X\times X \to Z$ be a function such that for some
$ \alpha >4$
$$N^{'}(\varphi_2(\frac{x}{2(3)},\frac{y}{2}),a)\geq
N^{'}(\varphi_2(\frac{x}{2(3)},y),\alpha a) $$ for all $x \in X,$
$y \in \{0,\frac{x}{3},\frac{4x}{3},\frac{-2x}{3},x\}$ and $a>0,$
and $\lim_{n \to \infty}N^{'}(4^n\varphi_2(2^{-n}x,2^{-n}y),a)=1$
for all $x, y \in X$ and $a
> 0.$ Let $f: X \to Y$ be
an even function with $f(0)=0$ satisfies $(2.2)$ for all $a
> 0$ and all $x, y \in X.$ Then there exists a unique quadratic mapping $Q:
X \to Y$ such that $$N(Q(x)-f(x),a)\geq
N^{''}_2(x,\frac{a(\alpha-4)}{6})$$ for all $x \in X$ and all $a
> 0,$ where
\begin{align*}
N^{''}_2(x,a):=min \{&N^{'}(\varphi_2(\frac{x}{3},\frac{x}{3}),a)
,N^{'}(\varphi_2(\frac{x}{3},x),a),\\
&N^{'}(\varphi_2(\frac{x}{3},\frac{4x}{3}),a)
,N^{'}(\varphi_2(\frac{x}{3},\frac{-2x}{3}),a)
,N^{'}(\varphi_2(\frac{x}{3},0),a)\}.
\end{align*}
\end{theorem}
\begin{proof}
The techniques are completely similar to those techniques  of
Theorem $2.1.$ Hence we present a sketch of proof. If we replace
$x$ by $\frac{x}{2^{n+1}}$ in $(2.11),$ then we have
$$N(4f(\frac{x}{2^{n+1}})-f(\frac{x}{2^n}),6a)\geq N^{''}_2(\frac{x}{2^{n+1}},a)$$
whence
$$N(4^{n+1}f(\frac{x}{2^{n+1}})-4^nf(\frac{x}{2^n}),6(4^n)a)\geq
N^{''}_2(\frac{x}{2^{n+1}},a)\eqno(2.19)$$ for all $x \in X$ and
$a>0.$ One can deduce
$$N(4^{n+m}f(\frac{x}{2^{n+m}})-4^mf(\frac{x}{2^m}),a)
\geq
N^{''}_2(x,\frac{a}{\sum_{i=1}^{n+m}\frac{6}{\alpha}(\frac{4}{\alpha})^i})\eqno(2.20)$$
for all $x \in X,n\geq0,m\geq0$ and $a>0.$ Hence, we conclude that
$\{4^nf(\frac{x}{2^n})\}$ is a Cauchy sequence in the fuzzy Banach
space $(Y,N).$ Therefore, there is a function $Q: X \to Y$ defined
by $Q(x):=N-\lim_{n \to \infty}4^nf(\frac{x}{2^n}).$ Employing
$(2.20)$ with $m=0$ we obtain $$N(Q(x)-f(x),a)\geq
N^{''}_2(x,\frac{a(4-\alpha)}{6})$$ for all $x \in X$  and all $a
> 0.$
\end{proof}
\begin{theorem}\label{t2}
Let $\varphi_3: X\times X \to Z$ be a function such that for some
$0 < \alpha < 2$
$$N^{'}(\varphi_3(2(\frac{x}{2}),2y))\geq N^{'}(\alpha
\varphi_3((\frac{x}{2}),y)) \eqno(2.21)$$ for all $x \in X,$ $y
\in \{x,\frac{x}{2},\frac{3x}{2},2x\}$ and $a>0,$ and $\lim_{n \to
\infty}N^{'}(\varphi_3(2^nx,2^ny),2^na)=1$ for all $x, y \in X$
and $a
> 0.$ Let $f: X \to Y$ be
an odd function  satisfying $(2.2)$ for all $a
> 0$ and all $x, y \in X.$ Then there exists a unique additive mapping $A:
X \to Y$ such that $$N(A(x)-f(x),a)\geq
N^{''}_3(x,\frac{a(2-\alpha)}{4})\eqno(2.22)$$ for all $x \in X$
and all $a> 0,$ where
\begin{align*}
N^{''}_3(x,a):= min \{ &N^{'}(\varphi_3(x,x),a)
,N^{'}(\varphi_3(\frac{x}{2}),a),\\
&N^{'}(\varphi_3(\frac{x}{2},2x),a)
,N^{'}(\varphi_3(\frac{x}{2},\frac{3x}{2}),a)\}.
\end{align*}
\end{theorem}
\begin{proof}
Replacing $y$ by $x$ in $(2.2),$ we get
$$N(f(3x)-3f(2x)+3f(x),a)\geq
N^{'}(\varphi_3(x,x),a)\eqno(2.23)$$ for all $x \in X$ and all
$a>0.$ Replacing $y$ by $3x$ in $(2.2),$ we get
$$N(f(5x)-f(4x)-f(2x)+f(x),a) \geq N^{'}(\varphi_3(x,3x),a)\eqno(2.24)$$
for all $x \in X$ and all $a>0.$ Putting $y=4x$ in $(2.2)$ we
obtain
$$N(f(6x)-f(5x)+f(3x)-3f(2x)+2f(x),a) \geq N^{'}(\varphi_3(x,4x),a)\eqno(2.25)$$
for all $x \in X$ and all $a>0.$ It follows from $(2.23),$
$(2.24)$ and $(2.25),$
\begin{align*}
N(f(6x)-f(4x)-f(2x),3a) \geq min \{ &N^{'}(\varphi_3(x,x),a)
,N^{'}(\varphi_3(x,3x),a),\\
&N^{'}(\varphi_3(x,4x),a)\} \hspace{5.cm}(2.26)
\end{align*}
for all $x \in X$ and all $a>0$. If we replace $x$ by
$\frac{x}{2}$ in $(2.26),$ then
\begin{align*}
N(f(3x)-f(2x)-f(x),3a)\geq
min \{ &N^{'}(\varphi_3(\frac{x}{2},\frac{x}{2}),a),N^{'}(\varphi_3(\frac{x}{2},2x),a),\\
&N^{'}(\varphi_3(\frac{x}{2},\frac{3x}{2}),a) \hspace{5.cm}(2.27)
\end{align*}
for all $x \in X$ and $a>0$. It follows from $(2.23)$ and
$(2.27),$
$$N(\frac{f(2x)}{2}-f(x),2a)\geq
N^{''}_3(x,a)\eqno(2.28)$$ for all $x \in X$ and all $a>0.$
Replacing $x$ by $2^nx$ in $(2.28),$ we get
$$N(\frac{f(2^{n+1}x)}{2}-f(2^nx),2a)\geq N^{''}_3(x,a)(2^nx,a)\eqno(2.29)$$ for all $x\in X$ and
 all $a>0.$ Using $(2.21)$ we get
$$N(\frac{f(2^{n+1}x)}{2}-f(2^nx),2a)\geq N^{''}_3(x,\frac{a}{\alpha^n})\eqno(2.30)$$
for all $x \in X$ and all $a>0.$ Replacing $a$ by $\alpha^na$ we
see that
$$N(\frac{f(2^{n+1}x)}{2^{n+1}}-\frac{f(2^nx)}{2^n},\frac{2a\alpha ^n}{2^n})
\geq N^{''}_3(x,a)\eqno(2.31)$$ for all $x \in X$ and all $a>0.$
It follows from
$\frac{f(2^nx)}{2^n}-f(x)=\sum_{i=0}^{n-1}\frac{f(2^{i+1}x)}{2^{i+1}}-\frac{f(2^ix)}{2^i}$
and $(2.31)$ that
$$N(\frac{f(2^nx)}{2^n}-f(x),\sum_{i=0}^{n-1}\frac{2a\alpha^i}{2^i})
\geq min
\bigcup_{i=0}^{n-1}\{N(\frac{f(2^{i+1}x)}{2^{i+1}}-\frac{f(2^ix)}{2^i},\frac{2a\alpha^i}{2^i})\}\\
\geq N^{''}_3(x,a)\eqno(2.32)$$  for all $x \in X$ and all $a>0.$
By replacing $x$ with $2^mx$ in $(2.32)$ we observe that
$$N(\frac{f(2^{n+m}x)}{2^{n+m}}-\frac{f(2^mx)}{2^m},\sum_{i=0}^{n-1}\frac{2a\alpha
^i}{2^{i+m}})\geq N^{''}_3(2^mx,a)\geq
N^{''}_3(x,\frac{a}{\alpha^m}),$$ whence
$$N(\frac{f(2^{n+m}x)}{2^{n+m}}-\frac{f(2^mx)}{2^m},\sum_{i=m}^{n+m-1}\frac{2a\alpha
^i}{2^{i}})\geq N^{''}_3(x,a)$$ for all $x \in X,$ $a>0$ and $m,n\geq 0.$\\
Hence
$$N(\frac{f(2^{n+m}x)}{2^{n+m}}-\frac{f(2^mx)}{2^m},a))
\geq N^{''}_3(x,\frac{a}{\sum_{i=m}^{n+m-1}\frac{2\alpha
^i}{2^{i}}})\eqno(2.33)$$ for all $x \in X,$  $a>0$ and $m,n\geq
0.$ Since $0<\alpha <2$ and
$\sum_{i=0}^{\infty}(\frac{\alpha}{2})^i< \infty$ the Cauchy
criterion for convergence and $(N_5)$ show that
$\{\frac{f(2^nx)}{2^n}\}$ is a Cauchy sequence in $(Y,N).$ Since
$(Y,N)$ is a fuzzy Banach space, this sequence converges to some
point $A(x)\in Y.$ So one can define the mapping $A:X \to Y$ by
$A(x):=N-\lim_{n \to \infty}\frac{f(2^nx)}{2^n}$ for all $x \in
X.$ Letting $m=0$ in $(2.33),$ we get
$$N(\frac{f(2^{n}x)}{2^{n}}-f(x),a)\geq
N^{''}_3(x,\frac{a}{\sum_{i=0}^{n-1}\frac{2\alpha
^i}{2^{i}}})\eqno(2.34)$$ for all $x \in X,$ and $a>0.$ Taking the
limit as $n \to \infty$ and using $(N_6)$ we get
$$N(A(x)-f(x),a)\geq N^{''}_3(x,a)$$
 for all $x \in X$ and all $a>0.$ Now, we show that $A$ is additive. Replace $x,y$ by
$2^nx,2^ny,$ respectively in $(2.2)$ to get
\begin{align*}
N(\frac{f(2^n(2x+y))}{2^n}&+\frac{f(2^n(2x-y))}{2^n}-\frac{f(2^n(x+y))}{2^n}
\\
&-\frac{f(2^n(x-y))}{2^n}-\frac{2f(2^n(2x))}{2^n}
+\frac{2f(2^nx)}{2^n},a)\\
&\geq N^{'}(\varphi_3(2^nx,2^ny),2^na)
\end{align*}
for all $x,y \in X$ and $a>0.$ Since $\lim_{n \to
\infty}N^{'}(\varphi_3(2^nx,2^ny),2^na)=1,$ then by Lemma $2.2$ of
\cite{Na} we get that the mapping $A: X \to Y$ is
additive.\\
To prove the uniqueness of $A,$ let $A^{'}:X \to Y$ be another
additive mapping satisfying $(2.22).$ Fix $x \in X.$ Clearly
$A(2^nx)=2^nA(x)$ and $A^{'}(2^nx)=2^nA^{'}(x)$ for all $n \in
\Bbb N.$ It follows from $(2.22)$ that
\begin{align*}
N(A(x)-A^{'}(x),a)&=N(\frac{A(2^nx)}{2^n}-\frac{A^{'}(2^nx)}{2^n},a)\\
&\geq min \{N(\frac{A(2^nx)}{2^n}-\frac
{f(2^nx)}{2^n},\frac{a}{2}),
N(\frac{f(2^nx)}{2^n}-\frac{A^{'}(2^nx)}{2^n},\frac{a}{2})\}\\
&\geq N^{''}_3(2^nx,\frac{2^na(2-\alpha)}{8}) \geq
N^{''}_3(x,\frac{2^na(2-\alpha)}{8\alpha^n})
\end{align*}
for all $x \in X$  and all $a>0.$\\
Since $\lim_{n \to
\infty}\frac{a(2^n)(2-\alpha)}{8\alpha^n}=\infty,$ we obtain
$\lim_{n \to
\infty}N^{''}_3(x,\frac{2^na(2-\alpha)}{8\alpha^n})=1.$ Therefore,
$N(A(x)-A^{'}(x),a)=1$ for all $a>0,$ whence $A(x)=A^{'}(x).$
\end{proof}
\begin{theorem}
Let $\varphi_4: X\times X \to Z$ be a function such that for some
$ \alpha >2$
$$N^{'}(\varphi_4(\frac{1}{2}(\frac{x}{2}),\frac{y}{2}),a)\geq N^{'}(
\varphi_4((\frac{x}{2}),y),\alpha a) $$ for all $x \in X,$ $y \in
\{x,\frac{x}{2},\frac{3x}{2},2x\}$ and $a>0,$ and $\lim_{n \to
\infty}N^{'}(2^n\varphi_4(2^{-n}x,2^{-n}y),a)=1$ for all $x, y \in
X$ and $a
> 0.$ Let $f: X \to Y$ be an odd function  satisfying $(2.2)$ for all $a
> 0$ and all $x, y \in X.$ Then there exists a unique additive mapping $A:
X \to Y$ such that $$N(A(x)-f(x),a)\geq
N^{''}_4(x,\frac{a(\alpha-2)}{4})$$ for all $x \in X$ and all $a>
0,$  where
\begin{align*}
N^{''}_4(x,a):= min \{ &N^{'}(\varphi_4(x,x),a)
,N^{'}(\varphi_4(\frac{x}{2}),a),\\
&N^{'}(\varphi_4(\frac{x}{2},2x),a)
,N^{'}(\varphi_4(\frac{x}{2},\frac{3x}{2}),a)\}.
\end{align*}
\end{theorem}
\begin{proof}
If we replace $x$ by $\frac{x}{2^{n+1}}$ in $(2.28),$ then we have
$$N(f(\frac{x}{2^n})-2f(\frac{x}{2^{n+1}}),a)\geq N^{''}_4(\frac{x}{2^{n+1}},a)$$
whence
$$N(2^nf(\frac{x}{2^n})-2^{n+1}f(\frac{x}{2^{n+1}}),(2^n)a)\geq
N^{''}_4(\frac{x}{2^{n+1}},a)$$ for all $x \in X$ and $a>0.$ One
can deduce
$$N(2^mf(\frac{x}{2^m})-2^{n+m}f(\frac{x}{2^{n+m}}),a)
\geq
N^{''}_4(x,\frac{a}{\sum_{i=1}^{n+m}\frac{1}{\alpha}(\frac{2}{\alpha})^i})\eqno(2.35)$$
for all $x \in X,n\geq0,m\geq0$ and $a>0.$ Hence, we conclude that
$\{2^nf(\frac{x}{2^n})\}$ is a Cauchy sequence in the fuzzy Banach
space $(Y,N).$ Therefore, there is a function $A: X \to Y$ defined
by $A(x):=N-\lim_{n \to \infty}2^nf(\frac{x}{2^n}).$ Employing
$(2.35)$ with $m=0$ we obtain $$N(A(x)-f(x),a)\geq N^{''}_4(x,a)$$
for all $x \in X$  and all $a
> 0.$\\ The rest of the proof is similar to the proof of theorem
$2.3.$
\end{proof}
We now prove our main theorem in paper.
\begin{theorem}\label{t2}
Let $\varphi: X\times X \to Z$ be a function such that for
some $0 < \alpha < 2$
$$N^{'}(\varphi(2(\frac{x}{2}),2y),a)\geq N^{'}(\alpha
\varphi(\frac{x}{2},y),a) $$ for all $x \in X,$ $y \in
\beta\{0,x,\frac{x}{2},\frac{4x}{3},\frac{-2x}{3},\frac{x}{3},\frac{3x}{2},2x\}$
and $a>0,$ and $\lim_{n \to
\infty}N^{'}(\varphi(2^nx,2^ny),2^na)=1$ for all $x, y \in X$ and
$a
> 0.$ Let $f: X \to Y$ with $f(0)=0$ be
a function  satisfying $(2.2)$ for all $a
> 0$ and all $x, y \in X.$ Then there exist a unique quadratic mapping $Q:
X \to Y$ and a unique additive mapping $A: X \to Y$  satisfying
$(1.3)$ and
$$N(Q(x)-A(x)-f(x),a)\geq N^{''}(x,a)\eqno(2.36)$$ for all $x \in X$ and all $a>
0,$ where $$N^{''}(x,a):= min
\{N^{''}_1(x,\frac{a(4-\alpha)}{12},N^{''}_3(x,\frac{a(2-\alpha)}{8}))\}$$
and $N^{''}_1,$ $N^{''}_3$ have been defined in Theorems $2.1$ and
$2.3,$ respectively.
\end{theorem}
\begin{proof}
Let $f_e(x)=\frac{f(x)+f(-x)}{2}$ for all $x \in X.$ Then
$f_e(0)=0,$ $f_e(-x)=f_e(x)$ and
\begin{align*}
N(f_e(2x&+y)+f_e(2x-y)-f_e(x+y)-f_e(x-y)-2f_e(2x)+2f_e(x),a)\\
&=N(\frac{1}{2}[f(2x+y)+f(2x-y)-f(x+y)-f(x-y)-2f(2x)+2f(x)]\\
&+\frac{1}{2}[f(-2x-y)+f(-2x+y)-f(-x-y)-f(-x+y)-2f(-2x)+2f(-x)],a)\\
&=N([f(2x+y)+f(2x-y)-f(x+y)-f(x-y)-2f(2x)+2f(x)]\\
&+[f(-2x-y)+f(-2x+y)-f(-x-y)-f(-x+y)-2f(-2x)+2f(-x)],2a)\\
&\geq min \{N^{'}(\varphi(x,y),a),N^{'}(\varphi(-x,-y),a)\}
\hspace{6cm}(2.37)
\end{align*}
for all $x,y \in X$ and $a>0.$  Hence, there exists a unique
quadratic function $Q: X \to Y$ satisfying
$$N(Q(x)-f_e(x),a)\geq N^{''}_1(x,\frac{a(4-\alpha)}{6}
\eqno(2.38)$$ for all $x \in X$  and all $a
> 0.$
Let $f_o(x)=\frac{f(x)-f(-x)}{2}$ for all $x \in X.$ Then
$f_o(0)=0,$ $f_o(-x)=-f_o(x)$ and
\begin{align*}
N(f_o(2x&+y)+f_o(2x-y)-f_o(x+y)-f_o(x-y)-2f_o(2x)+2f_o(x),a)\\
&\geq min \{N^{'}(\varphi(x,y),a),N^{'}(\varphi(-x,-y),a)\}
\end{align*}
for all $x,y \in X$ and $a>0.$ From Theorem $2.3,$ there exists a
unique additive function $A: X \to Y$ satisfying
$$N(A(x)-f_o(x),a)\geq N^{''}_3(x,\frac{a(2-\alpha)}{4})\eqno(2.39)$$
for all $x \in X$  and all $a> 0.$ Hence $(2.36)$ follows from
$(2.38)$ and $(2.39).$
\end{proof}

\end{document}